\documentclass[12pt]{amsart}
\usepackage[english]{babel}
\usepackage{amsmath}
\usepackage{amsthm}
\usepackage{amssymb}
\usepackage{mathrsfs}
\usepackage{enumerate}
\usepackage[notcite, final, notref]{showkeys}
\usepackage{dsfont}
\usepackage{tikz}
\usepackage[siunitx]{circuitikz}
\usetikzlibrary{arrows,calc,chains,shapes}
\newtheorem{theorem}{Theorem}[section]

\newtheorem{proposition}[theorem]{Proposition}
\newtheorem{observation}[theorem]{Observation}
\newtheorem{corollary}[theorem]{Corollary}

\theoremstyle{definition}
\newtheorem{remark}[theorem]{Remark}

\usepackage{xcolor}

\title[Tensor product of rational functions]
{realization of tensor-product and of\\[0.2cm]
Tensor-Factorization of rational functions}

\author[D. Alpay]{Daniel Alpay}
\address{(DA) 
Faculty of Mathematics, Physics, and Computation\\
Schmidt College of Science and Technology\\
Chapman University\\
One University Drive
Orange, California 92866\\
USA}
\email{alpay@chapman.edu}
\thanks{Daniel Alpay thanks the Foster G. and Mary McGaw Professorship in
Mathematical Sciences, which supported this research.}

\author[I. Lewkowicz]{Izchak Lewkowicz}
\address{(IL) Department of electrical engineering 
Ben-Gurion University of the Negev\\ P.O.B. 653\\ Beer-Sheva, 84105\\
Israel}
\email{izchak@bgu.ac.il}

\begin{document}

\begin{abstract}
We here first study the state space realization of a tensor-product of a pair of
rational functions.  At the expense of ``inflating" the dimensions, we
recover the classical expressions for realization of a regular product
of rational functions.
Then, under an additional assumption that the limit at infinity of a
given rational function exists and is equal to identity, an explicit formula for
a {\em tensor-factorization} of this function, is introduced.
\end{abstract}
\maketitle

\date{today}
\setcounter{equation}{0}
\section{Introduction}
\setcounter{equation}{0}
The problem of minimal factorization of matrix-valued rational functions of one complex variable has along history; see for instance
\cite{MR82e:47024,BarGohKaa1979,MR902611}. Less known seems to be
the counterpart of this problem when matrix product is replaced by tensor product.
More precisely, we study the following two problems: First, given two 
rational matrix-valued
functions $R_1$ and $R_2$ analytic at infinity, write a realization of the
tensor product $R_1\otimes R_2$ in terms of realizations of $R_1$ and $R_2$. Next, given
a matrix-valued rational function $R$ analytic at infinity, find its representations
as $R_1\otimes R_2$ where $R_1$ and $R_2$ are rational and analytic at infinity.\smallskip

To provide some motivation we note the following.
Tensor products play an important role in mathematics and quantum mechanics. In the latter case, a first example (see e.g. \cite[p. 162]{cohen31mecanique})
is the product of two wave functions, 
each belonging to a given Hilbert space, which belongs to the tensor product of the given Hilbert spaces; see e.g. \cite[Proposition 6.2, p. 111]{Neveu68} for the latter. 
Another example is the case of quantum states 
(positive matrices with trace equal to $1$; see e.g. \cite{MR2538183}). 
Given two states $M_1\in \mathbb C^{N_1\times N_1}$ and $M_2\in\mathbb C^{N_2\times N_2}$, of possibly different sizes, 
the tensor product $M_1\otimes M_2$ is still a state.
Note that if $M=M_1\otimes M_2$, one can recover $M_1$ and $M_2$ uniquely via the formula
\begin{equation}
d_1^*M_1c_1=\sum_{k=1}^{N_2}(d_1\otimes f_k)^*M(c_1\otimes f_k),\quad c_1,d_1\in\mathbb C^{N_1},
\label{inv}
\end{equation}
where $f_1,\ldots, f_{N_2}$ denotes an orthonormal basis for $\mathbb C^{N_2}$,
and similarly for $M_2$,
\begin{equation}
d_2^*M_2c_2=\sum_{k=1}^{N_1}(e_k\otimes d_2)^*M(e_k\otimes c_2),\quad c_2,d_2\in\mathbb C^{N_2},
\label{inv1}
\end{equation}
where now $e_1,\ldots, e_{N_1}$ is an orthonormal basis for $\mathbb C^{N_1}$.
See e.g. \cite[eq. (9.2.1) p. 97]{MR2538183}.\smallskip

If one starts from an arbitrary state $M\in\mathbb C^{N_1N_2\times N_1N_2}$ the 
matrices defined by \eqref{inv} and \eqref{inv1} 
will be states, called marginal states, but their tensor product
need not be equal to $M$.\smallskip

One can consider similar problems in the setting of functions. We focus the discussion 
on rational functions. If $R(z)$ is a
$\mathbb C^{N\times N}$-valued rational function and if $N=N_1N_2$, formulas
\eqref{inv} and \eqref{inv1} now define two rational functions  $R_A$ and $R_B$, 
respectively $\mathbb C^{N_1\times N_1}$ and $\mathbb C^{N_1\times N_1}$ -valued, via
\begin{equation}
\begin{split}
d_1^*R_A(z)c_1&=\sum_{k=1}^{N_2}(d_1\otimes f_k)^*R(z)(c_1\otimes f_k),\\
d_2^*R_B(z)c_2&=\sum_{k=1}^{N_2}(d_2\otimes f_k)^*R(z)(c_2\otimes f_k),\\
\end{split}
\end{equation}

If $R=R_1\otimes R_2$ where $R_1$ is $\mathbb C^{N_1\times N_1}$-valued and $R_2$ is $\mathbb C^{N_2\times N_2}$-valued, then these equations can be rewritten as
\begin{equation}
\begin{split}
R_A(z)&=R_1(z)\cdot({\rm Tr}\, R_2(z))\\
R_B(z)&=R_2(z)\cdot({\rm Tr}\, R_1(z))
\end{split}
\end{equation}
and so these equations basically solve the tensor factorization problem.\\

The purpose of this work is in a somewhat different direction; we would like to express
both tensor multiplication and tensor factorization of matrix-valued rational
functions using state space representations.\\

In the rest of this section we cite some known results.
Let $z_l$, $z_r$ (the subscript stands for ``left" and ``right") be a pair
of complex variables, and let $F_l(z_l)$, $F_r(z_r)$ be a pair of
$p_l\times m_l$, \mbox{$p_r\times m_r$-valued} rational functions,
respectively. Assume that neither has poles at infinity and denote by
$n_l$, $n_r$ the respective McMillan degrees. Thus, one can write
the rational functions and the respective realization as
\begin{equation}\label{eq:OriginalFunctions}
\begin{matrix}
F_l(z_l)=D_l+C_l(z_lI_{n_l}-A_l)^{-1}B_l
&&&
F_r(z_r)=D_r+C_r(z_rI_{n_r}-A_r)^{-1}B_r
\\~\\
R_{F_l}=
\begin{footnotesize}
\left(\begin{array}{c|c}
A_l&B_l\\
\hline
C_l&D_l\end{array}\right)
\end{footnotesize}
&&&
R_{F_r}=
\begin{footnotesize}
\left(\begin{array}{c|c}
A_r&B_r\\
\hline
C_r&D_r\end{array}\right).
\end{footnotesize}
\end{matrix}
\end{equation}
Recall that whenever $m_l=p_r$ the product $F_l(z_l)F_r(z_r)$ is
well-defined and its realization is given by\begin{footnote}{
Strictly speaking, in the references it was formulated for
$z_l=z_r=z$ i.e. for $F_l(z)F_r(z)$}\end{footnote} (see e.g.
\cite[Section 2.5]{BarGohKaaRan2010})
\begin{equation}\label{eq:RealizationOfProduct}
\begin{matrix}
R_{F_l F_r}
=
\begin{footnotesize}
\left(\begin{array}{cc|c}
A_l&B_lC_r&B_lD_r\\
0&A_r&B_r\\
\hline
C_l&D_lC_r&D_lD_r
\end{array}\right)
\end{footnotesize}
=
\begin{footnotesize}
\left(\begin{array}{c|c}
A&B\\
\hline
C&D
\end{array}\right)
\end{footnotesize}
=
\left(\begin{smallmatrix}
A_l&0&B_l\\
0&I_{n_r}&0\\
C_l&0&D_l
\end{smallmatrix}\right)
\left(\begin{smallmatrix}
I_{n_l}&0&0\\
0&A_r&B_r\\
0&C_r&D_r
\end{smallmatrix}\right),
\end{matrix}
\end{equation}
in the sense that
\begin{equation}
F_1(z_1)F_2(z_2)=D_1D_2+\begin{pmatrix} C_1&D_1C_2\end{pmatrix}\left(\begin{pmatrix}z_1I_{n_1}&0\\ 0& z_2I_{n_2}\end{pmatrix}-
\begin{pmatrix}A_1&0\\0&A_2\end{pmatrix}\right)^{-1}\begin{pmatrix}B_1D_2 \\ B_2\end{pmatrix}.
\end{equation}
If $z_l=z_r$ 
the sought
realization in \eqref{eq:RealizationOfProduct} is of McMillan degree
\[
n_l+n_r~.
\]
when minimal (roughly speaking when there is no pole-zero cancelation).
We next address ourselves to the {\em tensor product}\begin{footnote}
{In matrix theory circles known as the ``Kronecker product ", see
e.g. \cite[Section 4.2]{HornJohnson2}.}\end{footnote} of $F_l(z_l)$
and $F_r(z_r)$, resulting in $F_l\otimes F_r$, a ~\mbox{$p_lp_r\times
m_lm_r$-valued} rational function. Tensor product of rational
functions is discussed in \cite[Section 5.2]{FuhHel2015}.
\vskip 0.2cm

\noindent
So far for known results. In the next section we focus on $R_{F_l\otimes F_r}$,
the state space realization of $F_l\otimes F_r$. In Section 3 we set the framework for the main result, which is the factorization result presented in Section 4.

\section{Realization of a tensor-product of rational functions}

\noindent
We start with technicalities: We denote by boldface characters,
``inflated version" of the original ones, i.e.
\begin{equation}\label{eq:InflatedVersion1}
\begin{matrix}
\mathbf{A_l}:=A_l\otimes I_{p_r}
&&\mathbf{A_r}:=I_{m_l}\otimes A_r\\
\mathbf{B_l}:=B_l\otimes I_{p_r}
&&\mathbf{B_r}:=I_{m_l}\otimes B_r\\
\mathbf{C_l}:=C_l\otimes I_{p_r}
&&{\mathbf C_r}:=I_{m_l}\otimes C_r\\
\mathbf{D_l}:=D_l\otimes I_{p_r}
&&\mathbf{D_r}:=I_{m_l}\otimes D_r\\
\mathbf{F_l}(z_l):=\mathbf{C_l}\left(z_lI_{n_lp_r}-
\mathbf{A_l}\right)^{-1}\mathbf{B_l}+\mathbf{D_l}
&&
\mathbf{F_r}(z_r):=\mathbf{C_r}\left(z_lI_{m_ln_r}-
\mathbf{A_r}\right)^{-1}\mathbf{B_r}+\mathbf{D_r}~.
\end{matrix}
\end{equation}
We then show that at the expense of ``inflating" the dimensions one can replace a
{\em tensor product} by a {\em usual product}.

\begin{proposition}\label{Pn:RealizationTensor}
Let $F_l(z_r)$, $F_r(z_r)$ be a pair of $p_l\times m_l$, \mbox{$p_r\times m_r$-valued}
rational functions, of McMillan degree $n_l$, $n_r$, respectively, whose
realization is given in Eq. \eqref{eq:OriginalFunctions}. Following Eqs.
\eqref{eq:RealizationOfProduct} and \eqref{eq:InflatedVersion1},
one has that,
\begin{equation}\label{eq:Equivalence1}
R_{F_l\otimes F_r}=R_{\mathbf{F_l}\mathbf{F_r}}~.
\end{equation}
\end{proposition}

In order to go into details we shall repeatedly use the fact, see e.g.
\cite[Lemma 4.2.10]{HornJohnson2}, that for matrices $T\in{\mathbb C}^{n\times m}$, 
$X\in{\mathbb C}^{m\times l}$, $Y\in{\mathbb C}^{l\times p}$,
$Z\in{\mathbb C}^{p\times q}$
one has that
\begin{equation}\label{eq:Identity1}
TX\otimes YZ=(T\otimes Y)(X\otimes Z).
\end{equation}
We now explicitly compute the tensor product of $F_l(z_l)$ and
$F_r(z_r)$,
\[
\begin{smallmatrix}
F_l\otimes F_r&=&\left(D_l+C_l(z_lI_{n_l}-A_l)^{-1}B_l\right)\otimes
\left(D_r+C_r(z_rI_{n_r}-A_r)^{-1}B_r\right)
\\~\\~&=&
D_l\otimes D_r+D_l\otimes\left(C_r(z_rI_{n_r}-A_r)^{-1}B_r\right)
+\left(C_l(z_lI_{n_l}-A_l)^{-1}B_l\right)\otimes D_r
+\left(C_l(z_lI_{n_l}-A_l)^{-1}B_l\right)\otimes\left(
C_r(z_rI_{n_r}-A_r)^{-1}B_r\right) 
\end{smallmatrix}
\]
We next separately examine each block
\[
\begin{smallmatrix}
D_l\otimes\left(C_r(z_rI_{n_r}-A_r)^{-1}B_r\right)&=&
D_lI_{m_l}\otimes\left(C_r(z_rI_{n_r}-A_r)^{-1}B_r\right)
\\~&=&
\left(D_l\otimes\left(C_r(z_rI_{n_r}-A_r)^{-1}\right)\right)\left(
I_{m_l}\otimes B_r\right)
\\~&=&
\left(D_lI_{m_l}\otimes\left(C_r(z_rI_{n_r}-A_r)^{-1}\right)\right)\left(
I_{m_l}\otimes B_r\right)
\\~&=&
\left(D_l\otimes C_r\right)\left(I_{m_l}\otimes\left(
(z_rI_{n_r}-A_r)^{-1}\right)\right)
\left(I_{m_l}\otimes B_r\right)
\\~&=&
\left(D_l\otimes I_{p_r}\right)\left(I_{m_l}\otimes C_r\right)
\left(I_{m_l}\otimes\left(
(z_rI_{n_r}-A_r)^{-1}\right)\right)
\left(I_{m_l}\otimes B_r\right)
\\~&=&
\underbrace{
\left(D_l\otimes I_{p_r}\right)}_{\mathbf D_l}
\underbrace{\left(I_{m_l}\otimes C_r\right)}_{\mathbf C_r}
\left(z_rI_{m_ln_r}-\underbrace{I_{m_l}\otimes A_r}_{\mathbf{A_r}}\right)^{-1}
\underbrace{\left(I_{m_l}\otimes B_r\right)}_{\mathbf B_r}
\\~&=&
\mathbf{D_l}
\mathbf{C_r}
\left(z_rI_{m_ln_r}-
\mathbf{A_r}\right)^{-1}
\mathbf{B_r}
\end{smallmatrix}
\]
\[
\begin{smallmatrix}
\left(C_l(z_lI_{n_l}-A_l)^{-1}B_l\right)\otimes D_r&=&
\left(C_l(z_lI_{n_l}-A_l)^{-1}B_l\right)\otimes I_{p_r}D_r
\\~&=&
\left(C_l\otimes I_{p_r}\right)\left(\left((z_lI_{n_l}-A_l)^{-1}B_l
\right)\otimes D_r\right)
\\~&=&
\left(C_l\otimes I_{p_r}\right)\left(\left((z_lI_{n_l}-A_l)^{-1}B_l
\right)\otimes I_{p_r}D_r\right)
\\~&=&
\left(C_l\otimes I_{p_r}\right)\left(\left((z_lI_{n_l}-A_l)^{-1}\right
)\otimes I_{p_r}\right)
\left(B_l\otimes D_r\right)
\\~&=&
\left(C_l\otimes I_{p_r}\right)\left(\left((z_lI_{n_l}-A_l)^{-1}\right
)\otimes I_{p_r}\right)
\left(B_l\otimes I_{p_r}\right)\left(I_{m_l}\otimes D_r\right)
\\~&=&
\underbrace{\left(C_l\otimes I_{p_r}\right)}_{\mathbf C_l}\left((z_lI_{n_lp_r}-
\underbrace{A_l\otimes I_{p_r}}_{\mathbf A_l}\right)^{-1}
\underbrace{\left(B_l\otimes I_{p_r}\right)}_{\mathbf B_l}
\underbrace{\left(I_{m_l}\otimes D_r\right)}_{\mathbf D_r}
\\~&=&
\mathbf{C_l}\left((z_lI_{n_lp_r}-
\mathbf{A_l}\right)^{-1}\mathbf{B_l}\mathbf{D_r}
\end{smallmatrix}
\]
\[
\begin{smallmatrix}
\left(C_l(z_lI_{n_l}-A_l)^{-1}B_l\right)\otimes\left(C_r(z_rI_{n_r}-A_r)^{-1}
B_r\right) 
&=&
\left(C_l(z_lI_{n_l}-A_l)^{-1}B_lI_{m_l}\right)\otimes\left(
I_{p_r}C_r(z_rI_{n_r}-A_r)^{-1}B_r\right) 
\\~&=&
\left(C_l\otimes I_{p_r}\right)\left((z_lI_{n_l}-A_l)^{-1}B_l\right)\otimes
\left(C_r(z_rI_{n_r}-A_r)^{-1}\right)\left(
I_{m_l}\otimes B_r\right) 
\\~&=&
\left(C_l\otimes I_{p_r}\right)\left((z_lI_{n_l}-A_l)^{-1}B_lI_{m_l}\right)\otimes
\left(I_{p_r}C_r(z_rI_{n_r}-A_r)^{-1}\right)\left(
I_{m_l}\otimes B_r\right) 
\\~&=&
\left(C_l\otimes I_{p_r}\right)
\left((z_lI_{n_l}-A_l)^{-1}\otimes I_{p_r}\right)
\left(B_l\otimes C_r\right)\left(I_{m_l}\otimes\left((z_rI_{n_r}-A_r)^{-1}
\right)\right)
\left( I_{m_l}\otimes B_r\right) 
\\~&=&
\left(C_l\otimes I_{p_r}\right)
\left((z_lI_{n_lp_r}-
A_l\otimes I_{p_r}
\right)^{-1}
\left(B_l\otimes I_{p_r}\right)
\left(I_{m_l}\otimes C_r\right)
\left((z_rI_{m_ln_r}-
I_{m_l}\otimes A_r
\right)^{-1}
\left( I_{m_l}\otimes B_r\right)
\\~&=&
\mathbf{C_l}\left((z_lI_{n_lp_r}-\mathbf{A_l}\right)^{-1}
\mathbf{B_l}\mathbf{C_r}
\left((z_rI_{m_ln_r}-
\mathbf{A_r}\right)^{-1}
\mathbf{B_r}
\end{smallmatrix}
\]
Thus, one can write
\[
\begin{matrix}
F_l\otimes F_r&=
\underbrace{\scriptstyle D_l\otimes D_r}_{\mathbf{D}}
+
\left(\begin{smallmatrix}
\mathbf{C_l}&&\mathbf{D_l}\mathbf{C_r}
\end{smallmatrix}\right)
\left(\begin{smallmatrix}
\left((z_lI_{n_lp_r}-\mathbf{A_l}
\right)^{-1}
&
\left(z_lI_{n_lp_r}-\mathbf{A_l}
\right)^{-1}\mathbf{B_l}\mathbf{C_r}
\left(z_rI_{m_ln_r}-\mathbf{A_r}
\right)^{-1}
\\0&
\left(z_rI_{m_ln_r}-\mathbf{A_r}
\right)^{-1}
\end{smallmatrix}\right)
\left(\begin{smallmatrix}
\mathbf{B_l}\mathbf{D_r}\\~\\
\mathbf{B_r}
\end{smallmatrix}\right)
\\~\\~&=
\begin{smallmatrix}\mathbf{D}\end{smallmatrix}
+
\left(\begin{smallmatrix}
\mathbf{C_l}
&&\mathbf{D_l}\mathbf{C_r}
\end{smallmatrix}\right)
\left(~
\left(\begin{smallmatrix}
z_lI_{n_lp_r}&0\\
0&z_rI_{m_ln_r}
\end{smallmatrix}\right)
-
\left(\begin{smallmatrix}
\mathbf{A_l}
&&\mathbf{B_l}\mathbf{C_r}
\\
0&&\mathbf{A_r}
\end{smallmatrix}\right)~
\right)^{-1}
\left(\begin{smallmatrix}
\mathbf{B_l}\mathbf{D_r}\\~\\
\mathbf{B _r}
\end{smallmatrix}\right).
\end{matrix}
\]
Note that in particular
\[
D_l\otimes D_r=(D_lI_{m_l})\otimes(I_{p_r}D_r)=
\underbrace{(D_l\otimes I_{p_r})}_{\mathbf D_l}
\underbrace{(I_{m_l}\otimes D_r)}_{\mathbf D_r}=
\mathbf{D_l}\mathbf{D_r}={\mathbf D}.
\]
The realization of $F_l(z_l)\otimes F_r(z_r)$ can be compactly written as
\begin{equation}\label{eq:RealizationOfTensor}
R_{F_l\otimes F_r}
=
\begin{footnotesize}
\left(\begin{array}{cc|c}
\mathbf{A_l}
&
\mathbf{B_l}\mathbf{C_r}&\mathbf{B_l}\mathbf{D_r}\\
0&\mathbf{A_r}
&\mathbf{B_r}\\
\hline
\mathbf{C_l}&\mathbf{D_l}\mathbf{C_r}&\mathbf{D_l}\mathbf{D_r}
\end{array}\right)
\end{footnotesize}
=
\begin{footnotesize}
\left(\begin{array}{c|c}
\mathbf{A_o}&\mathbf{B_o}\\
\hline
\mathbf{C_o}&\mathbf{D}
\end{array}\right)
\end{footnotesize}
=
{\mathbf R},
\end{equation}
which is indeed in form of \eqref{eq:RealizationOfProduct},
\eqref{eq:InflatedVersion1}.
If $z_l=z_r$ and there is no pole-zero cancelation, the sought
realization in
\eqref{eq:RealizationOfTensor} is of McMillan degree
\[
n_lp_r+m_ln_r~.
\]
Note now that in a way similar to \eqref{eq:RealizationOfProduct},
one can factorize the realization in
\eqref{eq:RealizationOfTensor} as follows, 
\begin{equation}\label{eq:RealizationAsMatrixProduct1}
{\mathbf R}
=
\begin{footnotesize}
\left(\begin{array}{cc|c}
\mathbf{A_l}&
\mathbf{B_l}\mathbf{C_r}&\mathbf{B_l}\mathbf{D_r}\\
0&\mathbf{A_r}
&\mathbf{B_r}\\
\hline
\mathbf{C_l}&\mathbf{D_l}\mathbf{C_r}&
\mathbf{D}
\end{array}\right)
\end{footnotesize}
=
\left(\begin{smallmatrix}
\mathbf{A_l}&0&\mathbf{B_l}\\
0&I_{m_ln_r}&0\\
\mathbf{C_l}&0&\mathbf{D_l}
\end{smallmatrix}\right)
\left(\begin{smallmatrix}
I_{n_lp_r}&0&0\\
0&\mathbf{A_r}&\mathbf{B_r}\\
0&\mathbf{C_r}&\mathbf{D_r}
\end{smallmatrix}\right).
\end{equation}
We conclude this section by pointing out that Proposition 
\ref{Pn:RealizationTensor}
can be easily extended to more elaborate cases like
\[
F_a(z_a)\otimes F_b(z_b)\otimes F_c(z_c)\cdots
\]

\section{Realization of the inverse of a tensor product of rational functions}
\label{Sec:Inverse}

\noindent
For future reference, in this section we examine the realization
of the inverse of rational functions of the form
$F_l(z_l)\otimes F_r(z_r)$ studied in the previous section.
\vskip 0.2cm

\noindent
We first recall, see e.g. \cite[Theorem 2.4]{BarGohKaaRan2010}, in the
realization of the inverse a rational function: Namely if
\[
R_F=
\begin{footnotesize}
\left(\begin{array}{c|c}
A&B\\
\hline
C&D
\end{array}\right),
\end{footnotesize}
\]
is a realization of a square matrix-valued rational function $F(z)$,
then whenever $D$ is non-singular, $\left(F(z)\right)^{-1}$ is
well-defined almost everywhere, and a corresponding realization is given by,
\begin{equation}\label{eq:RealizaionOfInverse}
R_{F^{-1}}=
\begin{footnotesize}
\left(\begin{array}{c|c}
A^{\times}&B^{\times}\\
\hline
C^{\times}&D^{\times}
\end{array}\right)
\end{footnotesize}
=
\begin{footnotesize}
\left(\begin{array}{c|c}
A-BD^{-1}C&-BD^{-1}\\
\hline
D^{-1}C&D^{-1}
\end{array}\right).
\end{footnotesize}
\end{equation}
Next, whenever the above $F_l(z)$ and $F_r(z)$ are so that
\[
m_l=p_r\quad\quad
{\rm and}\quad\quad
p_l=m_r
\]
the product $F_l(z)F_r(z)$ is square, and whenever $D_lD_r$ is
non-singular\begin{footnote}
{this implies that $m_l=p_r\geq {\rm rank}(D_lD_r)=p_l=m_r$.}\end{footnote},
$\left(F_l(z)F_r(z)\right)^{-1}$ is well-defined almost everywhere,
and by combining \eqref{eq:RealizationOfProduct} together with 
\eqref{eq:RealizaionOfInverse} a corresponding realization is given by
\begin{equation}\label{eq:RealizationOfInverseProduct}
\begin{matrix}
R_{(F_l F_r)^{-1}}
=
\begin{footnotesize}
\left(\begin{array}{cc|c}
A_l^{\times}&0&B_l^{\times}\\
B_r^{\times}C_l^{\times}&A_r^{\times}&B_r^{\times}D_l^{-1}\\
\hline
D_r^{-1}C_l^{\times}&C_r^{\times}&D_r^{-1}D_l^{-1}
\end{array}\right)
\end{footnotesize}
=
\left(\begin{smallmatrix}
I_{n_l}&0&0\\
0&A_r^{\times}&B_r^{\times}\\
0&C_r^{\times}&D_r^{-1}
\end{smallmatrix}\right)
\left(\begin{smallmatrix}
A_l^{\times}&0&B_l^{\times}\\
0&I_{n_r}&0\\
C_l^{\times}&0&D_l^{-1}
\end{smallmatrix}\right).
\end{matrix}
\end{equation}
Similarly, whenever
\[
m_lm_r=p_lp_r,
\]
the rational function $F_l(z)\otimes F_r(z)$ is square and if 
$D_l\otimes D_r=\mathbf{D_l}\mathbf{D_r}=\mathbf{D}$
is non-singular, then $D_l$, $D_r$ are square, i.e.
\[
m_l=p_l\quad\quad\quad m_r=p_r
\]
and non-singular, see e.g. \cite[Theorem 4.2.15]{HornJohnson2}.
Thus, we shall denote hereafter by $m_l\times m_l$,
$m_r\times m_r$ the dimensions of $F_l$, $F_r$, respectively.
\vskip 0.2cm

\noindent
Under these conditions, the \mbox{$m_lm_r\times m_lm_r$-valued}
rational function, $\left(F_l(z)\otimes F_r(z)\right)^{-1}$
is almost everywhere defined.
\eqref{eq:RealizationOfTensor}, we next compute the realization of
$(F_l\otimes F_r)^{-1}$, 
\[
R_{(F_l\otimes F_r)^{-1}}
=
\left( \begin{smallmatrix}
\left(A_l\otimes I_{p_r}\right)-\left(B_l\otimes D_r\right)
\left(D_l\otimes D_r\right)^{-1}\left(C_l\otimes I_{p_r}\right)
&B_l\otimes C_r-\left(B_l\otimes D_r\right)\left(D_l\otimes D_r\right)^{-1}
\left(D_l\otimes C_r\right)
&-\left(B_l\otimes D_r\right)\left(D_l\otimes D_r\right)^{-1}\\
-\left(I_{m_l}\otimes B_r\right)\left(D_l\otimes D_r\right)^{-1}
\left(C_l\otimes I_{p_r}\right)
&
\left(I_{m_l}\otimes A_r\right)-\left(I_{m_l}\otimes B_r\right)
\left(D_l\otimes D_r\right)^{-1}\left(D_l\otimes C_r\right)
&
-\left(I_{m_l}\otimes B_r\right)\left(D_l\otimes D_r\right)^{-1}
\\
\left(D_l\otimes D_r\right)^{-1}\left(C_l\otimes I_{p_r}\right)
&
\left(D_l\otimes D_r\right)^{-1}\left(D_l\otimes C_r\right)
&
\left(D_l\otimes D_r\right)^{-1}
\end{smallmatrix}\right).
\]
Taking into account the fact that $D_l$ and $D_r$ are square
and non-singular,
the realization $R_{(F_l\otimes F_r)^{-1}}$ takes the form
\[
\begin{matrix}
R_{(F_l\otimes F_r)^{-1}}
&=&
\begin{footnotesize}
\left(\begin{array}{cc|c}
\left(A_l-B_lD_l^{-1}C_l\right)\otimes I_{p_r}
&0&
\left(-B_lD_l^{-1}\right)\otimes I_{p_r}\\
\left(D_l^{-1}C_l\right)\otimes\left(- B_rD_r^{-1}\right)
&
I_{m_l}\otimes\left(A_r-B_rD_r^{-1}C_r\right)
&
D_l^{-1}\otimes\left(- B_rD_r^{-1}\right)
\\
\hline
D_l^{-1}C_l\otimes D_r^{-1}
&
I_{m_l}\otimes D_r^{-1}C_r
&
\left(D_l\otimes D_r\right)^{-1}
\end{array}\right)\end{footnotesize}
\\~\\~&=&
\begin{footnotesize}
\left(\begin{array}{cc|c}
A_l^{\times}{\otimes} I_{p_r}
&0&
B_l^{\times}\otimes I_{p_r}\\
C_l^{\times}\otimes B_r^{\times}
&
I_{m_l}\otimes A_r^{\times}
&
D_l^{-1}\otimes B_r^{\times}
\\
\hline
C_l^{\times}\otimes D_r^{-1}
&
I_{m_l}\otimes C_r^{\times}
&
\left(D_l\otimes D_r\right)^{-1}
\end{array}\right)\end{footnotesize}
\\~\\~&=&
\begin{footnotesize}
\left(\begin{array}{cc|c}
A_l^{\times}\otimes I_{p_r}&0&\left(B_l^{\times}\otimes I_{p_r}\right)
\\
\left(I_{n_l}\otimes B_r^{\times}\right)\left(C_l^{\times}\otimes I_{p_r}\right)
&
I_{m_l}\otimes A_r^{\times}
&
\left(I_{m_l}\otimes B_r^{\times}\right)
\left(D_l^{-1}\otimes I_{n_r}\right)
\\
\hline
\left(I_{n_l}\otimes D_r^{-1}\right)
\left(C_l^{\times}\otimes I_{p_r}\right)
&
I_{m_lp_r}\otimes\left(I_{m_l}\otimes C_r^{\times}\right)
&
\left(I_{m_l}\otimes D_r^{-1}\right)
\left(D_l^{-1}\otimes I_{p_r}\right)
\end{array}\right)\end{footnotesize}
\\~\\~&=&
\begin{footnotesize}
\left(\begin{array}{cc|c}
\mathbf{A_l^{\times}}&0&\mathbf{B_l^{\times}}
\\
\mathbf{B_r^{\times}}\mathbf{C_l^{\times}}&
\mathbf{A_r^{\times}}&\mathbf{B_r^{\times}}\mathbf{D_l}^{-1}
\\
\hline
\mathbf{D_r^{-1}}\mathbf{C_l^{\times}}&\mathbf{C_r^{\times}}&
\mathbf{D_r^{-1}}
\mathbf{D_l^{-1}}
\end{array}\right)\end{footnotesize}
=
\begin{footnotesize}
\left(\begin{array}{c|c}
\mathbf{{A_o}^{\times}}&\mathbf{{B_o}^{\times}}\\
\hline
\mathbf{{C_o}^{\times}}&\mathbf{D^{\times}}
\end{array}\right)
\end{footnotesize}
=
\mathbf{R^{\times}},
\end{matrix}
\]
where the boldface entries are given by
\begin{equation}\label{eq:InflatedVersion2}
\begin{matrix}
\mathbf{A_l^{\times}}:=A_l^{\times}\otimes I_{p_r}
&\mathbf{A_r^{\times}}:=I_{m_l}\otimes A_r^{\times}\\
\mathbf{B_l^{\times}}:=B_l^{\times}\otimes I_{p_r}
&\mathbf{B_r^{\times}}:=I_{m_l}\otimes B_r^{\times}\\
\mathbf{C_l^{\times}}:=C_l^{\times}\otimes I_{p_r}
&\mathbf{C_r^{\times}}:=I_{m_l}\otimes C_r^{\times}\\
\mathbf{D_l^{-1}}=D_l^{-1}\otimes I_{p_r}
&\mathbf{D_r}^{-1}=I_{m_l}\otimes D_r^{-1}.
\end{matrix}
\end{equation}
One can conclude that
\[
R_{(F_l\otimes F_r)^{-1}}=R_{(\mathbf{F_l}\mathbf{F_r})^{-1}},
\]
and in a way similar to \eqref{eq:RealizationAsMatrixProduct1}, 
one can factorize the above realization as follows, 
\begin{equation}\label{eq:RealizationAsMatrixProduct2}
\begin{matrix}
\mathbf{R^{\times}}
=
\begin{footnotesize}
\left(\begin{array}{cc|c}
\mathbf{A_l^{\times}}&0&{\mathbf B_l^{\times}}
\\
\mathbf{B_r^{\times}}\mathbf{C_l^{\times}}&
\mathbf{A_r^{\times}}&\mathbf{B_r^{\times}}\mathbf{D_l^{-1}}
\\
\hline
\mathbf{D_r^{-1}}\mathbf{C_l^{\times}}&\mathbf{C_r^{\times}}&
\mathbf{D_r^{-1}}
{\mathbf D_l}^{-1}
\end{array}\right)\end{footnotesize}
=
\left(\begin{smallmatrix}
I_{n_lp_r}&0&0
\\
0&\mathbf{A_r^{\times}}&\mathbf{B_r^{\times}}
\\
0&\mathbf{C_r^{\times}}&\mathbf{D_r^{-1}}
\end{smallmatrix}\right)
\left(\begin{smallmatrix}
\mathbf{A_l^{\times}}
&0&
\mathbf{B_l^{\times}}
\\
0&I_{m_ln_r}&0
\\
\mathbf{C_l^{\times}}
&0&
\mathbf{D_l^{-1}}
\end{smallmatrix}\right).
\end{matrix}
\end{equation}

\section{Tensor-factorization of rational functions}
\setcounter{equation}{0}

We now address a more challenging question: Given $\mathbf{F}(z)$
and $\left(\mathbf{F}(z)\right)^{-1}$ (assuming that $\det\mathbf{F}(z)\not\equiv0$),
under what conditions and how, can it be ``tensor-factorized" to
{\em some} $F_l(z)$ and $F_r(z)$, namely the following
relation holds,
\begin{equation}\label{eq:TensorFactorization}
\mathbf{F}(z)=F_l(z)\otimes F_r(z).
\end{equation}
Note that here, we confine the
discussion to a single complex variable, i.e. \mbox{$z_l=z_r=z$}.
\vskip 0.2cm

\noindent
Note also that if \eqref{eq:TensorFactorization} holds, this is
true up to complex scaling i.e.,
\[
F_l(z)\otimes F_r(z)=c(z)F_l(z)\otimes{\scriptstyle\frac{1}{c(z)}}F_r(z)
\quad\quad\quad0\not=c(z)\in\mathbb{C}.
\]
We shall use this degree of freedom in the sequel.
\vskip 0.2cm

\noindent
We next recall in the following fact from matrix theory.
\vskip 0.2cm

\noindent
Let $\Pi_{\alpha}$, $\Pi_{\beta}$ be a pair of supporting projections of the
space $\mathbb{C}^{(\alpha+\beta)\times(\alpha+\beta)}$, i.e.
\begin{equation}\label{eq:SupporingProjections}
\begin{matrix}
\Pi_{\alpha}^2=\Pi_{\alpha}\\~\\ \Pi_{\beta}^2=\Pi_{\beta}
\end{matrix}
\quad\quad
\begin{matrix}
\Pi_{\alpha}\Pi_{\beta}=0_{\alpha+\beta}=\Pi_{\beta}\Pi_{\alpha}
\\~\\
\Pi_{\alpha}+\Pi_{\beta}=I_{\alpha+\beta}~.
\end{matrix}
\end{equation}
Such a pair of projections can be obtained by partitioning an arbitrary non-singular
$T\in\mathbb{C}^{(\alpha+\beta)\times(\alpha+\beta)}$ 
as follows.
\begin{equation}\label{eq:Projections}
\begin{matrix}
T^{-1}\left(\begin{smallmatrix}I_{\alpha}&0\\0&~~0_{\beta}\end{smallmatrix}\right)T
&:=&\Pi_{\alpha}
\\~\\
T^{-1}\left(\begin{smallmatrix}0_{\alpha}&0\\0&~~I_{\beta}\end{smallmatrix}\right)T
&=&
\Pi_{\beta}~.
\end{matrix}
\end{equation}
By using an isometry-like relation,
we next offer a simple way 
to ``deflate" matrix dimensions.

\begin{observation}\label{Ob:Factorization}
Given $M\in\mathbb{C}^{s\times q}$, denote
\[
\mathbf{M_l}:=M\otimes I_p
\quad\quad\quad
\mathbf{M_r}:=I_m\otimes M.
\]
For arbitrary $u\in\mathbb{C}^p$, $v\in\mathbb{C}^m$ normalized so that
$u^*u=1$, $v^*v=1$, one has that
\[
\begin{matrix}
\left(I_s\otimes u^*\right)\mathbf{M_l}
\left(I_q\otimes u\right)=M
\quad\quad{\rm and}\quad\quad
\left(v^*\otimes I_s\right)\mathbf{M_r}
\left(v\otimes I_q\right)=M.
\end{matrix}
\]
\end{observation}
\vskip 0.2cm

\noindent
Indeed, by twice applying \eqref{eq:Identity1} one obtains,
\[
\begin{matrix}
\left(I_s\otimes u^*\right)
\underbrace{\left(M\otimes I_p\right)}_{\mathbf{M_l}}
\left(I_q\otimes u\right)
&=&
\left(I_sMI_q\right)\otimes\underbrace{\left(u^*I_pu\right)}_{=1}
&=&M
\\
\left(v^*\otimes I_s\right)
\underbrace{\left(I_m\otimes M\right)}_{\mathbf{M_r}}
\left(v\otimes I_q\right)
&=&
\underbrace{\left(v^*I_m v\right)}_{=1}\otimes\left(I_s MI_q\right)
&=&M.
\end{matrix}
\]
We next apply the last observation to the variables here.

\begin{corollary}\label{Cy:Deflation}
For $u\in\mathbb{C}^{p_r}$, $v\in\mathbb{C}^{m_l}$, normalized so that
$u^*u=1$ and $v^*v=1$, the boldface characters in
\eqref{eq:InflatedVersion1} satisfy
\[
\begin{matrix}
A_l&=&(I_{n_l}\times u^*)\mathbf{A_l}(I_{n_l}\otimes u)
&&&
A_r&=&(v^*\otimes I_{n_l})\mathbf{A_r}(v\otimes I_{n_r})
\\
B_l&=&(I_{n_l}\times u^*)\mathbf{B_l}(I_{m_l}\otimes u)
&&&
B_r&=&(v^*\otimes I_{n_l})\mathbf{B_r}(v\otimes I_{m_r})
\\
C_l&=&(I_{p_l}\times u^*)\mathbf{C_l}(I_{n_l}\otimes u)
&&&
C_r&=&(v^*\otimes I_{p_l})\mathbf{C_r}(v\otimes I_{n_r})
\\
D_l&=&(I_{p_l}\times u^*)\mathbf{D_l}(I_{m_l}\otimes u)
&&&
D_r&=&(v^*\otimes I_{p_l})\mathbf{D_r}(v\otimes I_{m_r}).
\end{matrix}
\]
\end{corollary}

\noindent
We now return to the problem of ``tensor-factorization" in
\eqref{eq:TensorFactorization}.
We note that in place of $\mathbf{R}$ in
\eqref{eq:RealizationOfTensor}
and $\mathbf{R^{\times}}$ in 
\eqref{eq:RealizationAsMatrixProduct2}, the realization arrays associated with 
$\mathbf{F}$ and $\mathbf{F^{-1}}$,
are known only {\em up to a coordinate
transformation}, i.e. there exists, a non-singular matrix
$T\in{\mathbb C}^{(n_lp_r+m_ln_r)\times(n_lp_r+m_ln_r)}$
namely in \eqref{eq:SupporingProjections} and \eqref{eq:Projections}
\[
\alpha=n_lp_r\quad{\rm and}\quad\beta=m_ln_r~,
\]
so that the actual realization array is given by
\begin{equation}\label{eq:CoordinateTransformation1}
\left(\begin{smallmatrix}T&0\\0&I_{p_lm_r}\end{smallmatrix}\right)^{-1}
\mathbf{R}
\left(\begin{smallmatrix}T&0\\0&I_{p_lm_r}\end{smallmatrix}\right)
=
\begin{footnotesize}
\left(\begin{array}{c|c}
T^{-1}\mathbf{A_o}T&T^{-1}\mathbf{B_o}\\
\hline
\mathbf{C_o}T&\mathbf{D}
\end{array}\right)
\end{footnotesize}
=
\begin{footnotesize}
\left(\begin{array}{c|c}
\mathbf{A}&\mathbf{B}\\
\hline
\mathbf{C}&\mathbf{D}
\end{array}\right),
\end{footnotesize}
\end{equation}
and
\begin{equation}\label{eq:CoordinateTransformation2}
\left(\begin{smallmatrix}T&0\\0&I_{p_lm_r}\end{smallmatrix}\right)^{-1}
\mathbf{R^{\times}}
\left(\begin{smallmatrix}T&0\\0&I_{p_lm_r}\end{smallmatrix}\right)
=
\begin{footnotesize}
\left(\begin{array}{c|c}
T^{-1}{\mathbf{A_o}^{\times}}T&T^{-1}{\mathbf{B_o}^{\times}}\\
\hline
\mathbf{C_o^{\times}}T&\mathbf{D^{-1}}
\end{array}\right)
\end{footnotesize}
=
\begin{footnotesize}
\left(\begin{array}{c|c}
\mathbf{A^{\times}}&\mathbf{B^{\times}}\\
\hline
\mathbf{C^{\times}}&\mathbf{D^{-1}}
\end{array}\right).
\end{footnotesize}
\end{equation}
As in reality, the specific coordinate transformation, $T$ in 
\eqref{eq:CoordinateTransformation1} and
\eqref{eq:CoordinateTransformation2} is~ {\em unknown}~ one can conclude
that to extract $F_l(z)$ and $F_r(z)$ from \eqref{eq:TensorFactorization}
along with the realization arrays
in \eqref{eq:CoordinateTransformation1}, \eqref{eq:CoordinateTransformation2},
additional conditions are needed.

\begin{theorem}\label{Tm:Factorization}
Let ${\mathbf F}(z)$ be a given square matrix-valued rational function.
Assume that
\[
\lim\limits_{z~\longrightarrow~\infty}\mathbf{F}(z)=I.
\]
Let $~
\begin{footnotesize}
\left(\begin{array}{c|c}
\mathbf{A}&\mathbf{B}\\
\hline
\mathbf{C}&
I
\end{array}\right),
\end{footnotesize}~$
see \eqref{eq:CoordinateTransformation1}, and
$~
\begin{footnotesize}
\left(\begin{array}{c|c}
\mathbf{A^{\times}}&\mathbf{B^{\times}}\\
\hline
\mathbf{C^{\times}}&
I
\end{array}\right).
\end{footnotesize}
$
see \eqref{eq:CoordinateTransformation2}, be realizations of $\mathbf{F}(z)$
and of $\left(\mathbf{F}(z)\right)^{-1}$, respectively.
\vskip 0.2cm

\noindent
Substituting in \eqref{eq:SupporingProjections},
$\alpha=n_lm_r$ and $\beta=m_ln_r$,
assume also that there exists a pair of supporting projection to
$\mathbb{C}^{n_lm_r+m_ln_r}$
denoted by $\Pi_{n_lm_r}$ and $\Pi_{m_ln_r}$ so that
\begin{equation}\label{eq:SubspaceCondition}
\mathbf{A}\Pi_{n_lm_r}=\Pi_{n_lm_r}\mathbf{A}\Pi_{n_lm_r}
\quad\quad\quad
\mathbf{A^{\times}}\Pi_{m_ln_r}=\Pi_{m_ln_r}\mathbf{A^{\times}}
\Pi_{m_ln_r}~.
\end{equation}
Following the definition of the projections ${\Pi}_{n_lm_r}$ and
$~{\Pi}_{m_ln_r}$, see \eqref{eq:Projections} and \eqref{eq:SubspaceCondition},
along with Corollary \ref{Cy:Deflation}, for arbitrary
$u\in\mathbb{C}^{m_r}$, $v\in\mathbb{C}^{m_l}$, normalized so that
$u^*u=1$ and $v^*v=1$, we find it convenient to introduce the
following related projections\begin{footnote}{note that
$\hat{\Pi}_{n_lm_r}{\Pi}_{n_lm_r}=\hat{\Pi}_{n_lm_r}{\Pi}_{n_lm_r}=
\hat{\Pi}_{n_lm_r}$ and $~\hat{\Pi}_{m_ln_r}{\Pi}_{m_ln_r}=
{\Pi}_{m_ln_r}\hat{\Pi}_{m_ln_r}=\hat{\Pi}_{m_ln_r}$.}\end{footnote},
\begin{equation}\label{eq:HarProjection}
\hat{\Pi}_{n_lm_r}=
T^{-1}\left(\begin{smallmatrix}I_{n_l}\otimes uu^*&0\\
0&0_{m_ln_r}\end{smallmatrix}\right)T
\quad\quad\quad
\hat{\Pi}_{m_ln_r}=
T^{-1}\left(\begin{smallmatrix}0_{n_lm_r}&0\\
0&vv^*\otimes I_{m_l}\end{smallmatrix}\right)T
\end{equation}
Then, using \eqref{eq:RealizationOfTensor} and \eqref{eq:CoordinateTransformation1},
one can take in \eqref{eq:TensorFactorization} $\mathbf{F}=F_l\otimes F_r$ where,
\[
\begin{matrix}
F_l(z)&=&(I_{m_l}\otimes u^*)\mathbf{C}\hat{\Pi}_{n_lm_r}
\left(zI_{n_lm_r+m_ln_r}-\mathbf{A}\right)^{-1}
\hat{\Pi}_{n_lm_r}\mathbf{B}(I_{m_l}\otimes u)+I_{m_l}
\\~\\
F_r(z)&=&(v^*\otimes I_{m_r})\mathbf{C}\hat{\Pi}_{m_ln_r}
\left(zI_{n_lm_r+m_ln_r}-\mathbf{A}\right)^{-1}\hat{\Pi}_{m_ln_r}
\mathbf{B}(v\otimes I_{m_r})+I_{m_r}
\end{matrix}
\]
\end{theorem}
\vskip 0.2cm

\noindent
{\bf Proof :}~
First, recall (see Section \ref{Sec:Inverse}) that the assumption that
$D_l\otimes D_r=\mathbf{D_l}\mathbf{D_r}=\mathbf{D}$
is square non-singular, it implies that both $D_l$ and
$D_r$ are square non-singular. We shall thus denote the dimensions of
$F_l$ and $F_r$, by $m_l\times m_l$ and $m_r\times m_r$, respectively.
\vskip 0.2cm

\noindent
The assumption here that $\mathbf{D}=I_{m_lm_r}~$ implies (see e.e.
\cite[Theorem 4.2.12]{HornJohnson2}) that
\[
D_l=cI_{m_l}\quad\quad\quad D_r={\scriptstyle\frac{1}{c}}I_{m_r}
\quad\quad{\rm for~~some~~non-zero}\quad c\in\mathbb{C}.
\]
As already mentioned after \eqref{eq:TensorFactorization},
to simplify the exposition we shall take $c=1$.
\vskip 0.2cm

\noindent
Next, let $T$ in \eqref{eq:Projections}, \eqref{eq:CoordinateTransformation1},
\eqref{eq:CoordinateTransformation2} be the same so that the supporting
projections are $\Pi_{n_lm_r}$ and $\Pi_{m_ln_r}$. Next note that substituting 
\eqref{eq:RealizationOfTensor}, \eqref{eq:RealizationAsMatrixProduct2},
\eqref{eq:CoordinateTransformation1} and \eqref{eq:CoordinateTransformation2}
in condition \eqref{eq:SubspaceCondition} yields,
\[
\begin{matrix}
\mathbf{A}\Pi_{n_lm_r}=&
T^{-1}
\left(\begin{smallmatrix}\mathbf{A_l}&0\\0& 0_{m_ln_r}\end{smallmatrix}\right)
T&&&
\Pi_{m_ln_r}\mathbf{A}=&T^{-1}\left(\begin{smallmatrix}0_{n_lm_r}&
0\\0&\mathbf{A_r}\end{smallmatrix}\right)T
\\~\\
\mathbf{A^{\times}}\Pi_{m_ln_r}=&
T^{-1}
\left(\begin{smallmatrix}0_{n_lm_r}&0\\0&\mathbf{A_r^{\times}}
\end{smallmatrix}\right)T&&&
\Pi_{n_lm_r}\mathbf{A^{\times}}=&T^{-1}
\left(\begin{smallmatrix}\mathbf{A_l^{\times}}&0\\0& 0_{m_ln_r}
\end{smallmatrix}\right)T
\end{matrix}
\]
and thus in the sequel we shall use the two upper relations, i.e.
\[
\Pi_{n_lm_r}{\mathbf A}\Pi_{n_lm_r}=
T^{-1}
\left(\begin{smallmatrix}\mathbf{A_l}&0\\0& 0_{m_ln_r}\end{smallmatrix}\right)
T\quad\quad\quad
\Pi_{m_ln_r}\mathbf{A}\Pi_{m_ln_r}=
T^{-1}
\left(\begin{smallmatrix}0_{n_lm_r}&0\\0&\mathbf{A_r}\end{smallmatrix}\right)
T.
\]
We are now ready to recover $F_l(z)$, 
\[
\begin{matrix}
F_l(z)=C_l\left(zI_{n_l}-A_l\right)^{-1}B_l+I_{m_l}
\\~\\=
\underbrace{(I_{m_l}\otimes u^*)\mathbf{C_l}(I_{n_l}\otimes u)}_{C_l}
\underbrace{(I_{n_l}\otimes u^*)\left(zI_{n_lm_r}-\mathbf{A_l}\right)^{-1}
(I_{n_l}\otimes u)}_{\left(zI_{n_l}-A_l\right)^{-1}}
\underbrace{(I_{n_l}\otimes u^*)\mathbf{B_l}(I_{m_l}\otimes u)}_{B_l}
+I_{m_l}
\\~\\=
(I_{m_l}\otimes u^*)
\mathbf{C_l}
(I_{n_l}\otimes uu^*)
\left(zI_{n_lm_r}-\mathbf{A_l}\right)^{-1}
(I_{n_l}\otimes uu^*)
\mathbf{B_l}
(I_{m_l}\otimes u)+I_{m_l}
\\~\\=
(I_{m_l}\otimes u^*)
\mathbf{C_o}
\left(\begin{smallmatrix}I_{n_lm_r}\\~\\0_{m_ln_r\times n_lm_r}\end{smallmatrix}\right)
(I_{n_l}\otimes uu^*)
\left(\begin{smallmatrix}I_{n_lm_r}&0_{n_lm_r\times m_ln_r}\end{smallmatrix}\right)
\left(zI_{n_lm_r+m_ln_r}-\mathbf{A_o}\right)^{-1}
\\~\\
\times
\left(\begin{smallmatrix}I_{n_lm_r}\\~\\0\end{smallmatrix}\right)
(I_{n_l}\otimes uu^*)
\mathbf{B_l}(I_{m_l}\otimes u)+I_{m_l}
\\~\\
=
(I_{m_l}\otimes u^*)
\mathbf{C_o}
\left(\begin{smallmatrix}I_{n_l}\otimes uu^*&0\\
0&0_{m_ln_r}\end{smallmatrix}\right)
\left(zI_{n_lm_r+m_ln_r}-\mathbf{A_o}\right)^{-1}
\left(\begin{smallmatrix}I_{n_l}\otimes uu^*&0\\
0&0_{m_ln_r}\end{smallmatrix}\right)
\mathbf{B_o}
(I_{m_l}\otimes u)+I_{m_l}
\\~\\=
(I_{m_l}\otimes u^*)
\underbrace{
\mathbf{C_o}T}_{\mathbf C}\underbrace{
T^{-1}\left(\begin{smallmatrix}I_{n_l}\otimes uu^*&0\\
0&0_{m_ln_r}\end{smallmatrix}\right)T
}_{\hat{\Pi}_{n_lm_r}}\underbrace{T^{-1}
\left(zI_{n_lm_r+m_ln_r}-\mathbf{A_o}\right)^{-1}
T}_{\left(zI_{n_lm_r+m_ln_r}-\mathbf{A}\right)^{-1}}
\underbrace{T^{-1}
\left(\begin{smallmatrix}I_{n_l}\otimes uu^*&0\\
0&0_{m_ln_r}\end{smallmatrix}\right)T}_{\hat{\Pi}_{n_lm_r}}
\\~\\
\times
\underbrace{T^{-1}
\mathbf{B_o}}_{\mathbf B}
(I_{m_l}\otimes u)+I_{m_l}
\\~\\=
(I_{m_l}\otimes u^*)\mathbf{C}
\hat{\Pi}_{n_lm_r}
\left(zI_{n_lm_r+m_ln_r}-\mathbf{A}\right)^{-1}
\hat{\Pi}_{n_lm_r}
\mathbf{B}
(I_{m_l}\otimes u)+I_{m_l}~.
\end{matrix}
\]
Similarly, for $F_r(z)$
\[
\begin{matrix}
F_r(z)=C_r(zI_{n_r}-A_r)^{-1}B_r+I_{m_r}
\\~\\=
\underbrace{(v^*\otimes I_{m_r})\mathbf{C_r}(v\otimes I_{n_r})}_{C_r}
\underbrace{(v^*\otimes I_{n_l})(zI_{n_r}-\mathbf{A_r})^{-1}
(v\otimes I_{n_r})}_{(zI_{n_r}-A_r)^{-1}}
\underbrace{(v^*\otimes I_{n_l})\mathbf{B_r}(v\otimes I_{m_r})}_{B_r}
+I_{m_r}
\\~\\=
(v^*\otimes I_{m_r})\mathbf{C_r}(vv^*\otimes I_{m_l})
(zI_{n_r}-\mathbf{A_r})^{-1}(vv^*\otimes I_{m_l})
\mathbf{B_r}(v\otimes I_{m_r})+I_{m_r}
\\~\\=
(v^*\otimes I_{m_r})
\mathbf{C_o}
\left(\begin{smallmatrix}0_{n_lm_r\times m_ln_r}\\I_{m_ln_r}\end{smallmatrix}\right)
(vv^*\otimes I_{m_l})
\left(\begin{smallmatrix}0_{m_ln_r\times n_lm_r}&I_{m_ln_r}\end{smallmatrix}\right)
\left(zI_{n_lm_r+m_ln_r}-\mathbf{A_o}\right)^{-1}
\\~\\
\times
\left(\begin{smallmatrix}0\\~\\I_{m_ln_r}\end{smallmatrix}\right)
(vv^*\otimes I_{m_l})
\left(\begin{smallmatrix}0&I_{m_ln_r}\end{smallmatrix}\right)
\mathbf{B_o}(v\otimes I_{m_r})+I_{m_r}
\end{matrix}
\]
\[
\begin{matrix}
=(v^*\otimes I_{m_r})
\underbrace{\mathbf{C_o}T}_{\mathbf C}
\underbrace{
T^{-1}
\left(\begin{smallmatrix}0_{n_lm_r}&0\\0&vv^*\otimes I_{m_l}\end{smallmatrix}\right)
T
}_{\hat{\Pi}_{m_ln_r}}
\underbrace{T^{-1}
\left(zI_{n_lm_r+m_ln_r}-\mathbf{A_o}\right)^{-1}
T}_{\left(zI_{n_lm_r+m_ln_r}-\mathbf{A}\right)^{-1}}\underbrace{T^{-1}
\left(\begin{smallmatrix}0_{n_lm_r}&0\\0&vv^*\otimes I_{m_l}\end{smallmatrix}\right)
T}_{\hat{\Pi}_{m_ln_r}}
\\~\\
\times
\underbrace{T^{-1}\mathbf{B_o}}_{\mathbf B}(v\otimes I_{m_r})+I_{m_r}
\\~\\=
(v^*\otimes I_{m_r})\mathbf{C}\hat{\Pi}_{m_ln_r}
\left(zI_{n_lm_r+m_ln_r}-\mathbf{A}\right)^{-1}\hat{\Pi}_{m_ln_r}
\mathbf{B}(v\otimes I_{m_r})+I_{m_r}~.
\end{matrix}
\]
\qed
\vskip 0.2cm

\begin{remark}
At first sight, the assumptions in Theorem \ref{Tm:Factorization} seem
very restrictive. For persective recall that to factorize a given
rational function $F(z)$ to $F(z)=F_l(z)F_r(z)$, the assumptions are
virtually the
same\begin{footnote}{There they only assume $D$ is square non-singular,
but then only $F_l(z)D_r$ and $D_lF_r(z)$ are obtained.}\end{footnote},
see \cite[Section 2.5]{BarGohKaaRan2010}).
\end{remark}

\bibliographystyle{plain}
\def\cprime{$'$} \def\cprime{$'$} \def\cprime{$'$}
  \def\lfhook#1{\setbox0=\hbox{#1}{\ooalign{\hidewidth
  \lower1.5ex\hbox{'}\hidewidth\crcr\unhbox0}}} \def\cprime{$'$}
  \def\cprime{$'$} \def\cprime{$'$} \def\cprime{$'$} \def\cprime{$'$}
  \def\cprime{$'$}

\end{document}